\documentstyle[12pt,amssymb,epic,eepic]{amsart}

\def\draft{n}

\input psfig

\theoremstyle{plain}

\newtheorem{theorem}{Theorem}
\newtheorem{proposition}{Proposition}[section]
\newtheorem{lemma}[proposition]{Lemma}
\newtheorem{corollary}[proposition]{Corollary}

\theoremstyle{definition}
\newtheorem{definition}[proposition]{Definition}

\newtheorem{question}{Question}

\theoremstyle{remark}

\newtheorem{exercise}[proposition]{Exercise}

\newtheorem{remark}[proposition]{Remark}

\def\printname#1{
	\if\draft y
		\smash{\makebox[0pt]{\hspace{-0.5in}
			\raisebox{8pt}{\tt\tiny #1}}}
	\fi
}

\newcommand{\psdraw}[2]
        {\begin{array}{c} \hspace{-1.3mm}
	\raisebox{-4pt}{\psfig{figure=#1.ps,width=#2}}
	\hspace{-1.9mm}\end{array}}
\newlength{\standardunitlength}
\catcode`\@=11
\long\def\@makecaption#1#2{%
    \vskip 10pt
    \setbox\@tempboxa\hbox{
      \small\sf{\bfcaptionfont #1. }\ignorespaces #2}%
    \ifdim \wd\@tempboxa >\captionwidth {%
        \rightskip=\@captionmargin\leftskip=\@captionmargin
        \unhbox\@tempboxa\par}%
      \else
        \hbox to\hsize{\hfil\box\@tempboxa\hfil}%
    \fi}
\font\bfcaptionfont=cmssbx10 scaled \magstephalf
\newdimen\@captionmargin\@captionmargin=2\parindent
\newdimen\captionwidth\captionwidth=\hsize
\catcode`\@=12

\def\lbl#1{\label{#1}\printname{#1}}

\def\BZ{\Bbb Z}
\def\BQ{\Bbb Q}

\def\BC{\Bbb C}

\def\KK{\cal{K}}
\def\LL{\cal{L}}
\def\GG{\cal{G}}
\def\OO{\cal{O}}

\def\MM{\cal{M}}
\def\FF{\cal{F}}
\def\VV{\cal{V}}
\def\WW{\cal{W}}
\def\PP{\cal{P}}

\def\CMM{\cal{CM}}
\def\BMM{\cal{BM}}
\def\BBMM{\widetilde{\cal{BM}}}
\def\CCMM{\widetilde{\cal{CM}}}

\def\l{\lambda}

\def\as{algebraically split}
     
\def\CL{{\cal L}}
\def\CM{{\cal M}}

\def\FM#1{\cal F^{\text{}}_{#1}\cal M}
\def\FMM#1{\cal F^{\text{Ga}}_{#1}\cal M}
\def\FO#1{\cal F_{#1}\cal O}
\def\fti{finite type invariant}
\def\G#1{\cal G_{#1}}                         \def\GaroI{\cite{Ga}}
\def\GM#1{\cal G^{\text{}}_{#1}\cal M}        \def\GaroII{\cite{GL}}

\def\ihs{integral homology 3-sphere}          

                                              \def\OhII{\cite{Oh2}}
                                              \def\Rozansky{\cite{Rz1}}
\def\Q{\Bbb Q}                                
\def\uf{unit-framed}

\begin{document}


\title{On finite type 3-manifold invariants III: \newline
             Manifold Weight Systems}

\author{Stavros Garoufalidis}
\address{Department of Mathematics\\
        Massachusetts Institute of Technology\\
        Cambridge, MA 02139, U.S.A.}
\email{stavros@math.harvard.edu}
\thanks{The first named author was partially supported by NSF grant 
       DMS-95-05105.
This and related preprints can also be obtained 
by
       accessing the WEB in the address 
 {\tt 
http:\linebreak[0]//\linebreak[0]www.\linebreak[0]math.\linebreak
[
0]brown.\linebreak[0]edu/\linebreak[0]$\sim$stavrosg/} }

\author{Tomotada Ohtsuki}
\address{Department of Mathematical and Computing Sciences \\
Tokyo Institute of Technology \\          
Oh-okayama, Meguro-ku, Tokyo 152 \\          
Japan}
\email{tomotada@is.titech.ac.jp}

\date{This edition: March 12, 1997 \\ First edition: August 5, 1995 \\
      Fax number: 1 (617) 495 5132 \\
      Email: {\tt stavros@math.harvard.edu} \\
      Fax number: 81 (3) 5734 2714 \\
      Email: {\tt tomotada@is.titech.ac.jp}}

\maketitle

\begin{abstract}
The present paper is a continuation of \cite{Oh2} and \cite{GL} devoted to the
study of finite type invariants of \ihs s. We introduce the
notion of {\em manifold weight systems}, and show that  type $m$ 
invariants of \ihs s are determined (modulo invariants of type $m-1$) by 
their associated manifold weight systems. In particular we deduce a vanishing
theorem for finite type invariants. We show that the space of 
manifold weight systems forms a commutative, co-commutative Hopf algebra
and that the map from finite type invariants to manifold weight systems
is an algebra map. We conclude with better bounds for the graded space of
finite type invariants of \ihs s.
\end{abstract}

\tableofcontents


\section{Introduction}

\subsection{History}
\lbl{section.more.history}
The present paper is a continuation of \cite{Oh2} and \cite{GL} devoted to the
 study of finite type invariants of  oriented integral homology 3-spheres.

There are two main sources of motivation for the present work: (perturbative)
Chern-Simons theory in $3$ dimensions, and Vassiliev invariants of
knots in $S^3$.

Witten \cite{Wi} in his seminal paper, using path integrals (an infinite
dimensional ``integration'' method) introduced a topological quantum
field theory in $3$ dimensions whose Lagrangian was the Chern-Simons 
function on the space of all connections. The theory depends on the choice
of a semi-simple Lie group $G$ and an integer $k$. The expectation values
(in $\BC$) of the above mentioned quantum field theory yield invariants
of 3-manifolds (depending on $G$ and $k$) and invariants of knots in
3-manifolds (depending on $G$, $k$, and the choice of a representation
of $G$).

Though the above mentioned path integrals have not yet been defined, and an
analytic definition of the theory is not yet possible, shortly afterwards
many authors proposed (rigorous) combinatorial definitions 
for the above mentioned invariants of knots and 3-manifolds.
It turned out that the knot invariants (for knots in $S^3$) are values
at roots of unity for the various Jones-like polynomials,
what we call quantum invariants, of knots.

The results mentioned so far are non-perturbative, i.e. involve a path
integral over the space of all connections.

Perturbatively, i.e., in the limit $k \rightarrow \infty $ one expects
invariants of knots/3-manifolds that depend on the choice of a $G$ flat
connection. The $G$ flat connections are the critical points of the 
Chern-Simons function, and form a finite dimensional (however singular)
topological space. Due to the presence of a cubic term in the Chern-Simons 
function,  one expects contributions to the perturbative invariants
coming from trivalent graphs (otherwise called Feynman diagrams;
we will denote the set of Feynman diagrams by $\BMM$ below).
Such invariants have been defined by \cite{AxS1}, \cite{AxS2}.

Further we can expect that there should exist a universal quantum invariant
of 3-manifolds with values in the space of Feynman diagrams.

The second source of motivation comes from the theory of finite type
knot invariants.
{Finite type knot invariants} (or Vassiliev invariants)
were originally introduced by Vassiliev
\cite{Va}. For an excellent exposition including complete proofs of the
main theorems see \cite{B-N} and references therein.
We will proceed in analogy with features in Vassiliev invariants.

Finite type knot invariants have the following (rather appealing) features:

\begin{itemize}
\item 
     There is an {\em axiomatic} definition (over $\BZ$) (see \cite{B-N})
that resembles a ``difference formula'' of multivariable calculus.
They form a filtered commutative co-commutative algebra $\FF_{\star} \VV$.
\item  
     There is a map $\FF_{m} \VV \rightarrow \GG_m \WW $
to the space of {\em weight systems} $\GG_m \WW$. $\GG_m \WW$ is a 
{\em combinatorially} defined, finite dimensional vector space naturally
isomorphic to a few other vector spaces based on trivalent graphs 
(see \cite{B-N}).
\item  
     The above mentioned map is not one-to-one, however one has the
following short exact sequence: 
\begin{equation}\label{eq.V.exact}
0 \rightarrow \FF_{m-1} \VV \rightarrow \FF_m \VV \rightarrow \GG_m \WW
\rightarrow 0
\end{equation}
This is a general existence theorem for Vassiliev invariants due to
Kontsevich \cite{Ko}, which gives rise to the {\em universal} Vassiliev
invariant (with values in a space of chord diagrams).
\item  
     Examples of Vassiliev invariants are the derivatives (at $1$) of
quantum invariants of knots in $S^3$.
\item
     There are two contradictory conjectures: one that asserts that Vassiliev
invariants separate knots in $S^3$, and the other that asserts that all
Vassiliev invariants come from ``semi-simple Lie algebras''. There is favor
for each of them.
\end{itemize}

We observe similar features in the case of our finite type invariants
of 3-manifolds.
For the first feature, we define a commutative co-commutative
algebra structure in our space $\BMM$.
For the third feature, we have a short exact sequence 
in Theorem \ref{thm2} below,
though it might still be incomplete comparing to the above exact sequence.
For the fourth feature, we expect that
$\l_n$ defined in \OhII\ should be finite type, see Question \ref{que.2} below.

\subsection{Statement of the results}

We begin by introducing some notation and terminology which will be followed
in the rest of the paper.
Let $M$ be an oriented \ihs.
A {\it framed link} $\CL=(L,f)$ is an unoriented link 
$L=L_1 \cup L_2 \cup \cdots \cup L_n$ in $M$
with framing $f=(f_1,f_2,\cdots,f_n)$ (since $M$ is a \ihs, we can assume
that $f_i \in \BZ $).
A framed link $(L, f)$ is called {\em \uf} if $f_i =\pm1$ for all $i$.
A framed link $(L, f)$ is called {\em \as} if
the linking numbers between any two components of $L$ vanish,
{\em boundary} if each component of $L$ bounds a Seifert surface
such that the Seifert surfaces are disjoint from each other.

Let $\CM$ be the vector space over $\Q$ on the set of \ihs s.
The motivations from Chern-Simons theory and from finite type knot 
invariants described in Section \ref{section.more.history}
  gave rise to the following definition due to the second named
author of finite type invariants of \ihs:

\begin{definition}
\lbl{def.ftype}
A linear map $v : \CM \to \Q$ is a type $m$ invariant
in the sense of \OhII\ 
if it satisfies
$$
\sum_{\CL' \subset \CL} (-1)^{\# \CL'} v(M_{\CL'}) =0
$$
for every \ihs\ $M$ and
every \as\  \uf\ link $\CL$ of $m+1$ components in $M$,
where $\#\CL'$ denotes the the number of components of $\#\CL'$,
$M_{\CL'}$ the \ihs\ obtained by Dehn surgery on $M$ along $\CL'$
and the sum runs over all sublinks of $\CL$ including empty link.
Let $\FO m$ (resp. $ \OO $) denote the vector space of type $m$ 
(resp. finite type) invariants 
in the sense of \OhII\ . 
\end{definition}

\begin{remark}
It will be useful in the study of finite type invariants to introduce the
a decreasing filtration $\FF_{\star}^{}$ on  $\MM$.
For a framed link $\CL$ in $M$ we define $(M,\CL)$ by
\begin{equation}
\lbl{eq.notation}
(M,\CL) = \sum_{\CL' \subset \CL} (-1)^{\#\CL'} M_{\CL'} \in \CM.
\end{equation}
We can now define a decreasing filtration $\FF_{\star}^{}$ on the vector
space  $\MM$ as follows: $\FF_{m}^{} \MM$ is the subspace of $\MM$
spanned by $(M, \CL)$ 
for all \ihs\ $M$ and all \as\ 
\uf\ framed links of $n+1$ components.
Let $\GM n$ denote the associated graded vector space,
Note that $v \in \FO n$ if and only if $v(\FM{n+1})=0$.
This implies  that $\MM/\FF_{n+1} \MM$ is the dual vector space of
$\FO n$.
Note also that if $\LL \cup \KK$ denotes an \as\ \uf\ link (with a 
distinguished component of it, $\KK$) in a \ihs\ $M$,
then the following fundamental
relation holds:
\begin{equation}
\lbl{eq.fundamental}
(M, \LL \cup \KK)=(M, \LL)-(M_{\KK},\LL)
\end{equation}    
where we also denote by $\LL$ in $M_{\KK}$
the same framed link $\LL$ after doing Dehn surgery along $\KK$.
Note (trivially) that the left hand side of \eqref{eq.fundamental} 
involves links with one component more than the right hand side, an
observation that will be useful in the philosophical comment below. 
\end{remark}

\begin{remark}
A variation of  Definition \ref{def.basic} of finite type invariants of
\ihs s  was introduced by the first author \cite{Ga}: a type $m$ invariant
in that sense is a map $ v: \CM \to \Q$ such that $v(\FF_{m+1}^{Ga} \MM) =0$,
where $\FF_{m}^{Ga} \MM$ is the subspace of $\MM$ spanned by all pairs
$(M,\CL)$ for unit-framed boundary links $L$ in \ihs s $M$.
We will not study this definition in the present paper though.
\end{remark}

In analogy with the notion of Vassiliev invariants of knots, we introduce
the following notions:

\begin{definition}
\lbl{def.basic}
\begin{itemize}
\item 
     A {\em Chinese manifold character} (CMC) is a (possibly empty or
disconnected) graph
whose vertices are trivalent and oriented (i.e., one of the two
possible cyclic orderings of the edges emanating from such a vertex is
specified).
Let $\CMM$ denote the set of all Chinese manifold characters. $\CMM$
is a graded set (the degree of a Chinese manifold character is the number
of the edges of the graph).
\item
     An {\em extended Chinese manifold character} (ECMC) 
is a (possibly empty or disconnected) graph
whose vertices are either trivalent and oriented, or univalent. 
Let $\CCMM$ denote the set of all extended Chinese manifold characters. $\CCMM$
is a graded set with the same degree as the case of CMC.
\item
      Let
\begin{align}
\BMM &= span(\CMM)/ \{ AS, IHX \} \\
\BBMM &= span(\CCMM)/ \{ AS, IHX, IS \}
\end{align}
Here $AS$ and $IHX$ are the relations referred to in Figures \ref{figure.AS}
and \ref{figure.IHX}, and $IS$ is the set of extended Chinese manifold 
characters containing a component who is an interval, i.e., a graph with
two vertices and a single edge.
$\BMM$ and $\BBMM$ inherits a grading from $\CCMM$ and $\CMM$ respectively.
We denote by $\iota$ a linear map of $\BMM$ to $\BBMM$
which is induced by the inclusion of $\CMM$ to $\CCMM$.
\item
     A {\em manifold weight system} of degree $m$ is a map $W: \GG_m \BMM 
\rightarrow \BQ$. The set of manifold weight systems of degree $m$ 
is denoted by $\GG_m \WW \MM$.
\end{itemize}
\end{definition}

With the above notation, let us recall the following map introduced by the
second named author in \cite{Oh2}:
\begin{equation}
\lbl{eq.graphM}
\tilde{O}_{m}^{\star} : \GG_m \CCMM \rightarrow \GG_m \MM
\end{equation}
is defined as follows:
for a Chinese manifold character $\Gamma$ with $m$ edges, we consider the 
ribbon graph obtained by replacing every 
trivalent vertex with an oriented small disk,
and every edge by an oriented band.
If an edge has a trivalent vertex in its end,
the corresponding band is attached to the corresponding disk
preserving the orientation as in Figures
\ref{f2} and 
\ref{figure.halftwist}.
The result is an \as\ link $L(\Gamma)$ in $S^3$. We define
$\tilde{O}_{m}^{\star} (\Gamma)$ to be the image of 
$(L(\Gamma), (1,1, \dots, 1 ))
\in \FF_m \MM$ under the projection map  
$\FF_m \MM \rightarrow \GG_m \MM$.
Note that it is non-trivial to show that the map (\ref{eq.graphM})
is well defined.

\begin{figure}[htpb]
$$ \psdraw{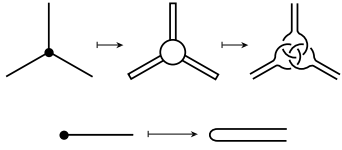}{1.8in}
$$
\caption{From vertex-oriented graphs to links in $S^3$.}\lbl{f2}
\end{figure}

\begin{figure}[htpb]
$$ \psdraw{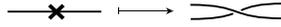}{1.5in}
$$
\caption{A marking of an edge and the part of the corresponding 
link.}\label{figure.halftwist}
\end{figure}

With these preliminaries, we recall  the following theorem due to the second
author:
\begin{theorem}\cite{Oh2} 
\lbl{thm.basic}
     The map \eqref{eq.graphM} is well defined and {onto}.
\end{theorem}

We can now state the main result of this paper:

\begin{theorem}
\footnote{The theorem has recently been improved by
the work of T.T.Q. Le, 
see the   Appendix. }
\lbl{thm2}
\begin{itemize}
\item
     The composition of the map \eqref{eq.graphM} with the {deframing map}
$F$ in \eqref{eq.frame} of Definition \ref{def.frame} 
 descends to a well defined map
\begin{equation}
\lbl{eq.graphD}
\iota \circ \tilde{O}_m^\star : \GG_m \BBMM \rightarrow \GG_m \MM
\end{equation}
We denote $\iota \circ \tilde{O}_m^\star $ by $O_m^\star$.
\item
     We have a surjection:
\begin{equation}\label{eq.Omstar}
 O_m^\star : \GG_m \BMM \rightarrow \GG_m \MM
\end{equation}
\item
     Dually, we get a map
$ O_m :
\FF_m \OO \rightarrow \GG_m \WW \MM$.
This map is not one-to-one, however one has the following exact sequence:
\begin{equation}\label{eq.O.exact}
0 \rightarrow \FF_{m-1} \OO \rightarrow \FF_m \OO \rightarrow \GG_m \WW \MM
\end{equation}
i.e.,
3-manifold invariants are determined in terms of their associated
manifold  weight systems.
\end{itemize}
\end{theorem}

Using the following lemma, 
which will be proved in Section \ref{section.vanish}
\begin{lemma}
\lbl{lemma.finite}
$\GG_m \CCMM/\{AS,IS \}$ (and therefore, $\GG_{m} \BBMM$ too) is a zero
dimensional vector space  if $m$ is not a multiple of $3$.
  Furthermore, 
$\GG_{3m} \BBMM$ is generated by Chinese manifold characters each
connected component of which is either the $Y$ component\footnote{
a $Y$ component is a graph with 1 vertex and 3 edges,
as in the letter $Y$}, or a trivalent graph
with no univalent vertices.
\end{lemma} 

and Theorem \ref{thm2} we obtain the following corollary:
\begin{corollary}
\lbl{corollary.vanish}
If $m$ is not a multiple of $3$, then $
\GG_m \OO=0$. In any case, it follows that $\GG_m \OO$ is a finite dimensional
vector space. 
\end{corollary}

\begin{remark}
The above Corollary \ref{corollary.vanish} was also observed by
\cite{GL}. It follows by the $AS$ relation (which itself follows by 
Theorem 4.1 of \cite{Oh2}
(see also \cite{GL})). In other words, the proof of Lemma \ref{lemma.finite} 
and Corollary \ref{corollary.vanish} and does not need the
extra $IHX$ relation  of Theorem \ref{thm2}.
\end{remark}

There is more structure on the vector spaces  $ \OO$ and $\BMM$
 which we now describe.
Using the pointwise multiplication of 3-manifold invariants, we observed
in \cite{Ga}
that there is a map $ \FF_m \OO \otimes \FF_n \OO \rightarrow \FF_{m+n}
\OO$, giving $ \OO$ the structure of a (filtered) commutative
algebra.

\begin{proposition}
\lbl{proposition.algebra}
\begin{itemize}
\item
     $\BMM$ (and therefore, $ \WW \MM$ as well) 
has a naturally defined multiplication $\cdot$ and commultiplication
$\Delta$, which together
make it a commutative, co-commutative Hopf algebra. 
By the structure theorem of Hopf algebras, (see \cite{Sw})
 it follows that $\BMM$ is 
the symmetric algebra on the (graded) set of the primitive elements 
\begin{equation}
\PP(\BMM) =\{ a \in \BMM :
\Delta(a)=a \otimes 1 + 1 \otimes a \}
\end{equation} in
$\BMM$.
\item
     The set of primitive elements $\PP(\BMM)$
is the set of
{\em connected} Chinese manifold characters.
\item
     Furthermore, the map
$ O_m : \FF_m \OO \rightarrow \GG_m \WW \MM $ is an algebra map.
\end{itemize}
\end{proposition}

\begin{corollary}
\lbl{corollary.calculations}\footnote{
The corollary has recently been improved by
the work of T.T.Q. Le, see the  Appendix  }
We have the following evaluation of 
dimensions of the graded vector spaces in low degrees.
$$
\vbox{
\offinterlineskip
\halign{\strut#&&\vrule#&\quad\hfil#\hfil\quad\cr
\noalign{\hrule}
&& $n$  && $0$ && $3$ && $6$ && $9$ && $12$     &\cr
\noalign{\hrule}
&& $dim \GG_n \OO$  && 
$1$ && $1$ && $1\le\cdot\le2$ && $1\le\cdot\le3$ && $1\le\cdot\le5$ &\cr
\noalign{\hrule}
&& $dim \GG_n \BMM$ && $1$ && $1$ && $2$ && $3$ && $5$  &\cr
\noalign{\hrule}
&& $dim \GG_n \PP(\BMM)$ && $1$ && $1$ && $1$ && $1$ && $2$  &\cr 
\noalign{\hrule}
}}
$$
\end{corollary}

\begin{pf}[of Corollary \ref{corollary.calculations}]
The third and fourth lines in the table follow by a direct calculation.
A computer version of the above calculation appears in \cite{B-N2}. The
lower bounds on the second line follow from the fact that $\GG_3 \OO$ is
a one dimensional vector space (spanned by the Casson invariant, as shown
in \cite{Oh2}) and the fact that $\GG_\star \OO$ is a commutative algebra.
The upper bounds follow from the third line and Theorem \ref{thm2}.
\end{pf}

The first author conjectured in \GaroI\ that $\FMM n=\FM{3n}$;
this implies a claim that
$\GM n=0$ for $n$ not a multiple of $3$,
which was further discussed by Rozansky in \Rozansky.
Corollary \ref{corollary.vanish} gives the positive answer to the claim.

\subsection{Plan of the proof}

The present paper consists of two parts, Sections \ref{section.more.proof}
and \ref{section.combi}.

Section \ref{section.more.proof} concentrates with 3-dimensional topology.
Using a restricted set of Kirby moves for surgery presentations of
\ihs s, we show Theorem \ref{thm2}. For the convenience of the reader,
we divide the proof in subsections \ref{sub.AS} and \ref{sub.IHX}.
We refer the reader to Section \ref{section.philo} for a philosophical 
comment on the
proof of Theorem \ref{thm2}. In Section \ref{section.vanish} we prove
Lemma \ref{lemma.finite}, and thus deduce Corollary \ref{corollary.vanish}.

Section \ref{section.combi} concentrates on the combinatorial aspects
of trivalent graphs. We show Proposition \ref{proposition.algebra} and
deduce Corollary \ref{corollary.calculations}.

Finally, in Section \ref{section.que} we pose some questions related to
finite type invariants of \ihs s.

\subsection{A philosophical comment on the proof of Theorem \ref{thm2}}
\lbl{section.philo}
In case the proof of Theorem \ref{thm2} is not too clear, the reader may
find useful the following philosophical comment. In order to state it,
let us introduce the following terminology: we call a {\em blow up} (of an
element in $\FF_\star \MM$), the move of replacing the right hand side 
of the fundamental relation (\ref{eq.fundamental}) by the left hand side
in an expression of the above mentioned element. Similarly, we call a
{\em blow down} the opposite move. Note that a blow up (respectively, a 
blow down) increases (respectively, decreases) the
number of the components of the link by one. With this terminology
we can restate remark $3.3$ of \cite{GL} as follows: surgical equivalence
is the relation generated by a sequence $B_1^+ B_1^- B_2^+ B_2^-
B_3^+ B_3^- \cdots$ (where $B_i^+$ are blow ups, $B_i^-$ are blow downs).
A mnemonic way for convincing oneself about that, is keeping track of
the number of components of the links in the various proofs involved.
Similarly, the relation in Theorem $4.1$ of \cite{Oh2} (for a precise
statement, see Theorem $5$ of \cite{GL}) and the $AS$ relation  of the
present paper is proven using a sequence
$B_1^- B_1^+ B_2^- B_2^+ B_3^- B_3^+ \cdots$. 
The $IHX$ relation in the present paper is proven using a sequence
$B_1^- B_2^- B_1^+ B_1^+ \cdots$.
Of course, blow ups do not commute with blow downs, and as a  result of
this non-commutativity we obtain the $IHX$ relation.  

\subsection{Questions}
\lbl{section.que}

\begin{question}\footnote{
The question has recently being answered positively  by
T.T.Q. Le, see the  Appendix. }
\lbl{que.1}
Is it true that the map 
$ O_m : \FF_m \OO \rightarrow \GG_m \WW \MM $ is onto, i.e., does every
manifold weight system integrate to a \ihs\ invariant?
\end{question}

\begin{remark}
Note that the analogous statement of Question \ref{que.1} for knot invariants 
is true, but harder than the previous statements about weight systems.
Note also that a positive answer to Question \ref{que.1} implies that 
$ \OO$ is
a commutative co-commutative Hopf algebra (with commultiplication defined 
by the connected sum of \ihs), and therefore a symmetric
algebra in a graded set of generators.
\end{remark}

\begin{question}
\lbl{que.2}
In the case of $sl_2$, is $\l_n$ defined in \OhII\ finite type?
\end{question}

\begin{question}
Is it true that finite type invariants of \ihs s separate them?
\end{question}

\begin{remark}
The above two questions may well be contradicting each other, as in the
case of knot invariants.
\end{remark}

Much remains to be done. 

\subsection{Acknowledgment}
The authors wish to thank the Mathematical Institute at Aarhus for the 
warm hospitality during which the last part of the paper was written.
In particular, they wish to thank J. Andersen for organizing an
excellent conference on finite type 3-manifold invariants, and for bringing
together participants from all over the world. 
Especially they wish to thank P. Melvin and the anonymous referee
 for pointing out to us
an  insufficient explanation in the proof of Lemma \ref{lemma.more}
 and  the {\tt Internet}
for supporting  numerous electronic  communications.


\section{3-dimensional topology}
\lbl{section.more.proof}

In this section we prove Theorem \ref{thm2}.
Our proof exploits the fact that diffeomorphic \ihs\ can be represented
in different ways as Dehn surgery on framed links in $S^3$. Though we do
not know of a complete set of moves that relate two surgery presentations
(within the category of \ihs s), we can still prove Theorem \ref{thm2}.

 We will show that the map
$\tilde{O}_m^\star: \GG_m \CCMM \rightarrow \GG_m \MM$, after {\em deframing}
 factors through a map 
$ \GG_m \BMM/ \{AS, IHX, FR \} \rightarrow \GG_m \MM$.

We begin with the following remark on drawing figures.
\begin{remark}
\lbl{rem.figures}
In all figures, the parts of the graphs and links not shown are
assumed identical. We call a figure {\em homogenous} if
all graphs (or links) shown have the same number of components;
otherwise we call it {\em inhomogenous}.
Given a  figure, the graphs or links drawn in it
represent elements in the graded space $\MM/\FF_{n+1}\MM$ under the map
of equation \eqref{eq.graphM}, where $n$ is the maximum of the number
of components of the links shown on the figure. Notice that in a homogenous
figure of $n$-component graphs or links, it is obvious that each element
is a well defined element of $\MM/\FF_{n+1}\MM$ under the map
of equation \eqref{eq.graphM}, whereas in a non-homogenous figure
it needs to be shown.
Some examples of homogenous Figures are \ref{figure.mark}, \ref{figure.Tbreak},
\ref{figure.2half}, \ref{figure.break}, \ref{f8}, \ref{f9}, \ref{figure.AS},
and \ref{f11} shown below. 
Some examples of inhomogenous Figures are \ref{f12} and  \ref{f21}
shown below.
\end{remark}

\subsection{Removing the marking}
\lbl{sub.remove}
We begin with the following lemma:

\begin{lemma}
\lbl{lemma.mark}
If two Chinese manifold characters with $m$ edges 
differ by their marking as shown in
Figure \ref{figure.mark}, then they have the same image in $\GG_m \MM $
under the map $\tilde O^{\star}_m$. 
\end{lemma}

\begin{figure}[htpb]
$$ \psdraw{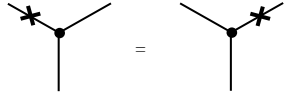}{1.5in}
$$
\caption{Moving the  marking on adjacent edges of a graph.}\lbl{figure.mark}
\end{figure}

\begin{pf}
We denote $L_1$, $L_2$ and $L_3$ be the three components
of a framed link which are images of three edges around the trivalent vertex
through the map $\tilde O^{\star}_m$.
In the same argument in \OhII,
we can regard $L_3$ as an element $[x_1,x_2]$ in the fundamental group
of the complement of $L_1 \cup L_2$
which is generated by two meridians $x_1$ and $x_2$.
If we make a mark on the first (resp. second) edge near the trivalent vertex,
the element becomes $[x_1^{-1},x_2]$ (resp. $[x_1,x_2^{-1}]$).
Since these two elements are conjugate in the fundamental group,
they express the same element in $\GM{m}$ as in \OhII.
This implies the relation in Figure \ref{figure.mark}.
\end{pf}

\begin{lemma}
\lbl{lemma.break}    
We have a relation in $\G{\star}\CCMM$ shown in Figure \ref{figure.Tbreak}.
\end{lemma}

\begin{figure}[htpb]
$$ \psdraw{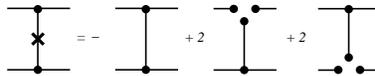}{2.0in}
$$
\caption{Removing a marking of a graph.}\lbl{figure.Tbreak}
\end{figure}

\begin{pf}
This relation is a direct conclusion of a relation
in $\GM{\star}$ shown in Figure \ref{figure.2half},
which can be obtained in the same way as in \OhII. It also follows
from Theorem $5$ of \cite{GL}.
\end{pf}

\begin{figure}[htpb]
$$ \psdraw{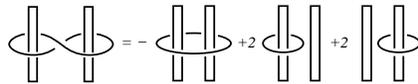}{2.2in}
$$
\caption{Untying  a half twist}\lbl{figure.2half}
\end{figure}

\subsection{A graphical representation of the deframing map
$F: \BBMM \to \BBMM$}
\lbl{sub.deframe}

In order to give a cleaner form of Lemmas \ref{lemma.mark} and 
\ref{lemma.break} (and also being motivated by Theorem $5$ of \cite{GL}),
we introduce white vertices, whose graphical definition
is shown in Figure \ref{figure.break}.
More precisely, we can give an alternative definition
through the following deframing map.

\begin{definition}
\lbl{def.frame}
The {\em deframing map}
\begin{equation}
\lbl{eq.frame}
F :\CCMM \rightarrow  \CCMM
\end{equation}
is defined as follows: for a manifold Chinese character $\Gamma$, let
\begin{equation}
F(\Gamma)=\sum\begin{Sb} c \subseteq v_3(\Gamma) \end{Sb}
(-1)^{|c|} \Gamma^c
\end{equation}
where the summation is over the set of all subsets of the set
$v_3(\Gamma)$ of trivalent vertices of $\Gamma$, and $\Gamma^c$ is the graph 
obtained
by splitting each trivalent vertex in $c$ with $3$ univalent ones, as in
Figure \ref{figure.break}.
 Here $|c|$ stands for the cardinality of the set $c$. 
\end{definition} 

We now compose the map \eqref{eq.graphM} with the
{deframing map} \eqref{eq.frame} of Definition \ref{def.frame}.

The deframing map $F$ is a map between two copies of $\CCMM$.
In order to distinguish these two copies,
we use the following conventions;
we draw  an extended Chinese manifold character graph
which belongs to the first $\CCMM$ (source of $F$)
by a graph with white trivalent vertices ( $\circ$ ),
whereas an extended Chinese manifold character graph
which belongs to the second $\CCMM$ (image of $F$)
is drawn   by a graph with black trivalent vertices ( $\bullet$ ).
By identifying these two copies of $\CCMM$ through $F$, 
we obtain the relation between $\circ$ and
$\bullet$ vertices as shown in Figure \ref{figure.break}.

\begin{figure}[htpb]
$$ \psdraw{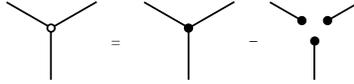}{1.9in}
$$
\caption{The definition of a $\circ$ vertex.}\lbl{figure.break}
\end{figure}

With respect to this substitution,
Lemmas \ref{lemma.mark} and \ref{lemma.break} become
the following two lemmas.

\begin{lemma}
\lbl{lemma.newmark}
A mark can move beyond a white trivalent vertex.
\end{lemma}

\begin{pf}
This lemma is immediately obtained from Lemma \ref{lemma.mark}
 by the definition
of white vertex, noting that a mark near a univalent vertex can be removed.
\end{pf}

\begin{lemma}
\lbl{lemma.newbreak} 
The relation in Figure \ref{f8} holds.
\end{lemma}

\begin{figure}[htpb]
$$ \psdraw{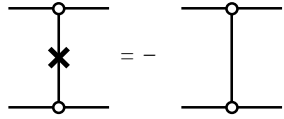}{1.5in}
$$
\caption{An equivalent form of Lemma \ref{lemma.break}.}\lbl{f8}
\end{figure}

\begin{pf}
This lemma is obtained from Lemma \ref{lemma.break}
by the definition of white vertex and
the fact that a graph including a connected component 
with one edge and two univalent vertices is equivalent to zero.
\end{pf}

\subsection{A vanishing lemma}
\lbl{section.vanish}

\begin{lemma} 
\lbl{lemma.UNI}
If a graph has a connected component
containing a univalent vertex,
no black vertices and at least two white vertices,
then it is equal to zero in $\GG_\star \MM$.
\end{lemma}

\begin{pf}
This lemma follows from the calculation in Figure \ref{f9}
using Lemmas \ref{lemma.newmark} and \ref{lemma.newbreak}.
\end{pf}

\begin{figure}[htpb]
$$ \psdraw{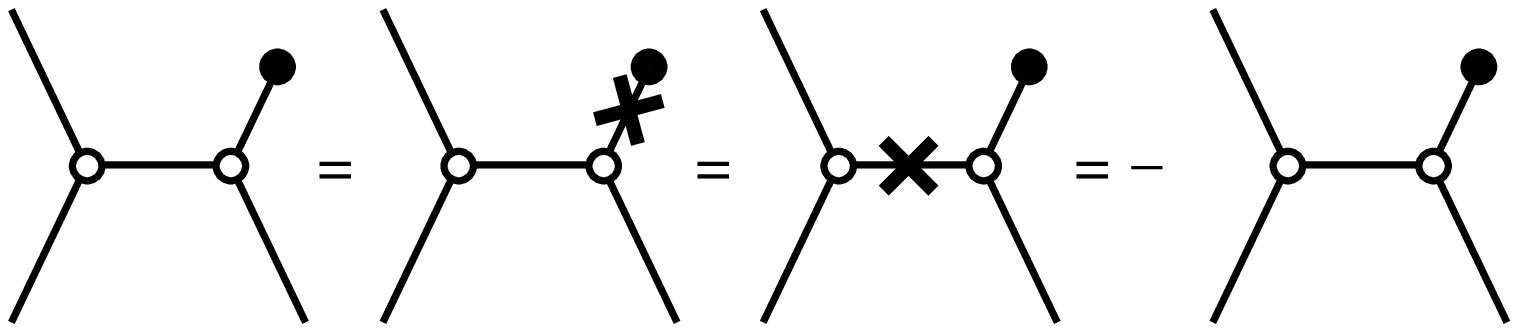}{2.2in}
$$
\caption{Proof of Lemma \ref{lemma.UNI}}\lbl{f9}
\end{figure}

We can now give a proof of Lemma \ref{lemma.finite} as follows:

\begin{pf}[of Lemma \ref{lemma.finite}]
The space $\GG_n \CCMM$ is spanned by graphs with
white trivalent vertices and univalent vertices.
Let $\Gamma$ be a such graph.

If a graph $\Gamma$ contains a univalent vertex, then 
it is equivalent to zero unless the univalent vertex belongs to
a $Y$ component.
Hence we can assume every connected component of $\Gamma$ is 
either a $Y$ graph, or a 
trivalent graph with white vertices and no univalent vertices.
Since the number of edges in any trivalent graph is divisible by $3$,
we obtain this lemma.
\end{pf}

\subsection{Proof of the  $AS$ relation}
\lbl{sub.AS}

\begin{proposition}
\lbl{prop.AS} 
The relation in Fig 2.9 holds,
which is called {\em AS (anti-\linebreak
symmetry)} relation.
Here we use \lq\lq blackboard cyclic order'',
that is, we assume  that each vertex has clockwise cyclic order
when it is depicted in a plane.
\end{proposition}

\begin{figure}[htpb]
$$ \psdraw{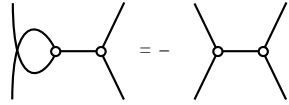}{1.5in}
$$
\caption{The $AS$ relation with an extra vertex}\lbl{figure.AS}
\end{figure}

\begin{pf}
We show a proof using Lemmas \ref{lemma.newmark} and \ref{lemma.newbreak}
 in Figure \ref{f11}.
\end{pf}

\begin{figure}[htpb]
$$ \psdraw{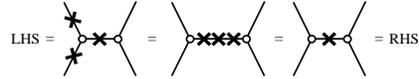}{2.3in}
$$
\caption{Proof of Proposition \ref{prop.AS}}\lbl{f11}
\end{figure}

\begin{remark}
Up to now, with the terminology of \cite{GL}, we only used the relations
of surgical equivalence and the fundamental relation of Theorem
$4.1$ of \cite{Oh2} (for a precise reformulation, see Theorem $5$
of \cite{GL}) in order to
show the $AS$ relation of Figure \ref{figure.AS}.
\end{remark}

\begin{remark}
One word of caution though: in the 3-manifold graphs, the $Y$ graph does not
vanish, whereas in the Vassiliev invariants graphs, the $Y$ graph vanishes.
The reason is that the 3-manifold $AS$ relation needs an external vertex,
otherwise it is a {\em symmetry} relation. This was observed in \cite{GL}
too, and seems to be responsible for the existence of the Casson invariant.
\end{remark}

\subsection{Proof of the $IHX$ relation}
\lbl{sub.IHX}
   
In order to prove $IHX$ relation, we begin with  the following lemma.

\begin{lemma}
\lbl{lemma.preparation} 
The relation in Figure \ref{f12} holds,
where by a band we mean either the  empty set or two parallel strings
with opposite orientations in a component of a framed link 
as shown in Figure \ref{f13}.
\end{lemma}

\begin{figure}[htpb]
$$ \psdraw{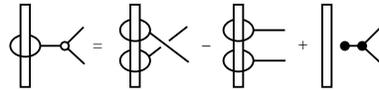}{2.0in}
$$
\caption{Breaking a white vertex. The 
uni-trivalent graphs in the first and fourth (resp. second
and third) parts of the figure have $n$ (resp. $n-1$) components.
The figure represents an identity in $\MM/\FF_{n+1} \MM$}\lbl{f12}
\end{figure}

\begin{figure}[htpb]
$$ \psdraw{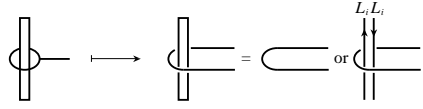}{2.2in}
$$
\caption{The definition of a band}\lbl{f13}
\end{figure}

\begin{pf}[Proof of Lemma \ref{lemma.preparation}]
Consider two framed $n$-component links shown in Figure \ref{f14}, 
whose framings are all $+1$.
By $(\#1)$ and $(\#2)$ we mean elements in $\GM{n}$,
where we express $(S^3,\CL) \in \GM{\star}$ by a picture of $\CL$
according to remark \ref{rem.figures}.

\begin{figure}[htpb]
$$ \psdraw{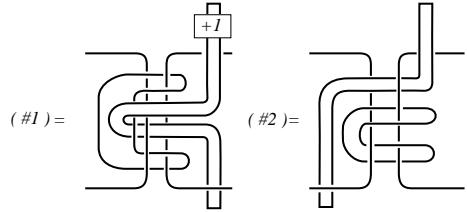}{2.5in}
$$
\caption{The definition of $(\#1)$ and $(\#2)$}\lbl{f14}
\end{figure}

Let $\CL_1$ (resp. $\CL'_1$) be the middle component of the framed link 
$\CL$ (resp. $\CL'$) which expresses $(\#1)$ (resp. $(\#2)$).
Note that we can obtain $\CL$ 
by taking handle slide of the band in $\CL'$ along $\CL'_1$,
which means $(S^3_{\CL_1},\CL-\CL_1)=(S^3_{\CL'_1},\CL'-\CL'_1)$.
Since 
\begin{eqnarray*}
(S^3,\CL) & = & (S^3,\CL-\CL_1)-(S^3_{\CL_1},\CL-\CL_1) \\
(S^3,\CL') & = & (S^3,\CL'-\CL_1')-(S^3_{\CL_1'},\CL'-\CL_1')
\end{eqnarray*}
we can calculate $(\#1)-(\#2)$ by eliminating the middle components
as in Figure \ref{f15}. 

\begin{figure}[htpb]
$$ \psdraw{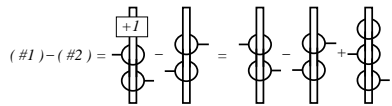}{2.1in}
$$
\caption{The calculation of $(\#1)-(\#2)$. An inhomogenous
 figure that represents
an identity in $\MM/\FF_{n+1}\MM$.}\lbl{f15}
\end{figure}

On the other hand we can deform $(\#1)$ and $(\#2)$ into graphs respectively.
We show a calculation for $(\#1)$ in Figure \ref{f16}
where we use Lemma \ref{lemma.more.prep} below to obtain the last term.
In a similar way we can obtain the corresponding graph to $(\#2)$
as shown in Figure \ref{f17}.

\begin{figure}[htpb]
$$ \psdraw{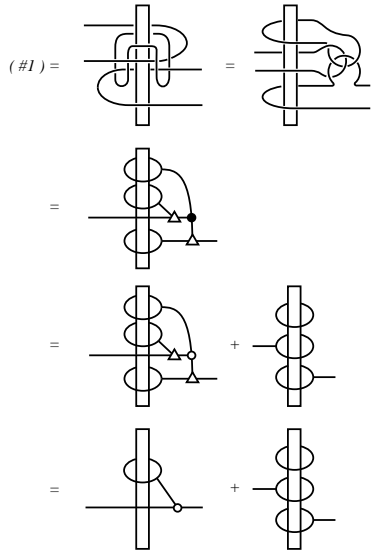}{2in}
$$
\caption{The calculation of $(\#1)$}\lbl{f16}
\end{figure}

\begin{figure}[htpb]
$$ \psdraw{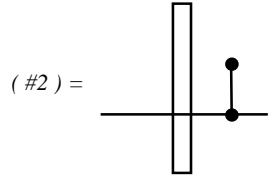}{1.5in}
$$
\caption{The calculation of $(\#2)$}\lbl{f17}
\end{figure}

Combining the results of
 Figures \ref{f15}, \ref{f16} and \ref{f17} we obtain the 
required formula.
\end{pf}

\begin{lemma}
\lbl{lemma.more.prep}  
The relation in Figure \ref{f18} holds. 
\end{lemma}

\begin{figure}[htpb]
$$ \psdraw{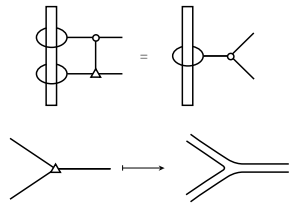}{1.5in}
$$
\caption{The definition of a triangle (shown on the lower part of the figure) 
and an identity
in $\GG_n\MM$ among white and triangle vertices of $n$-component links
(shown on the upper part of the figure).}\lbl{f18}
\end{figure}

\begin{pf}
Since we can change an order of two winding parts
in a framed link of $n$ components in $\GM n$,
we obtain the formula in Figure \ref{f19}.

\begin{figure}[htpb]
$$ \psdraw{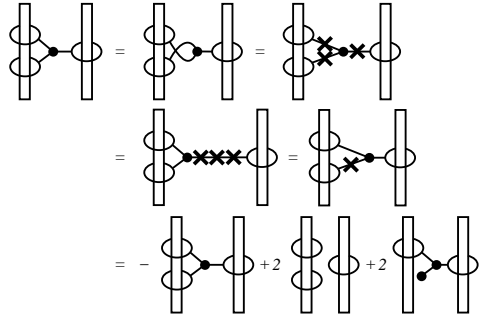}{2.5in}
$$
\caption{A vertex between two bands}\lbl{f19}
\end{figure}

Using Lemma 3.4 in \GaroII\ and the above formula,
we have the calculation in Figure \ref{f20}, which shows the required formula.
\end{pf}

\begin{figure}[htpb]
$$ \psdraw{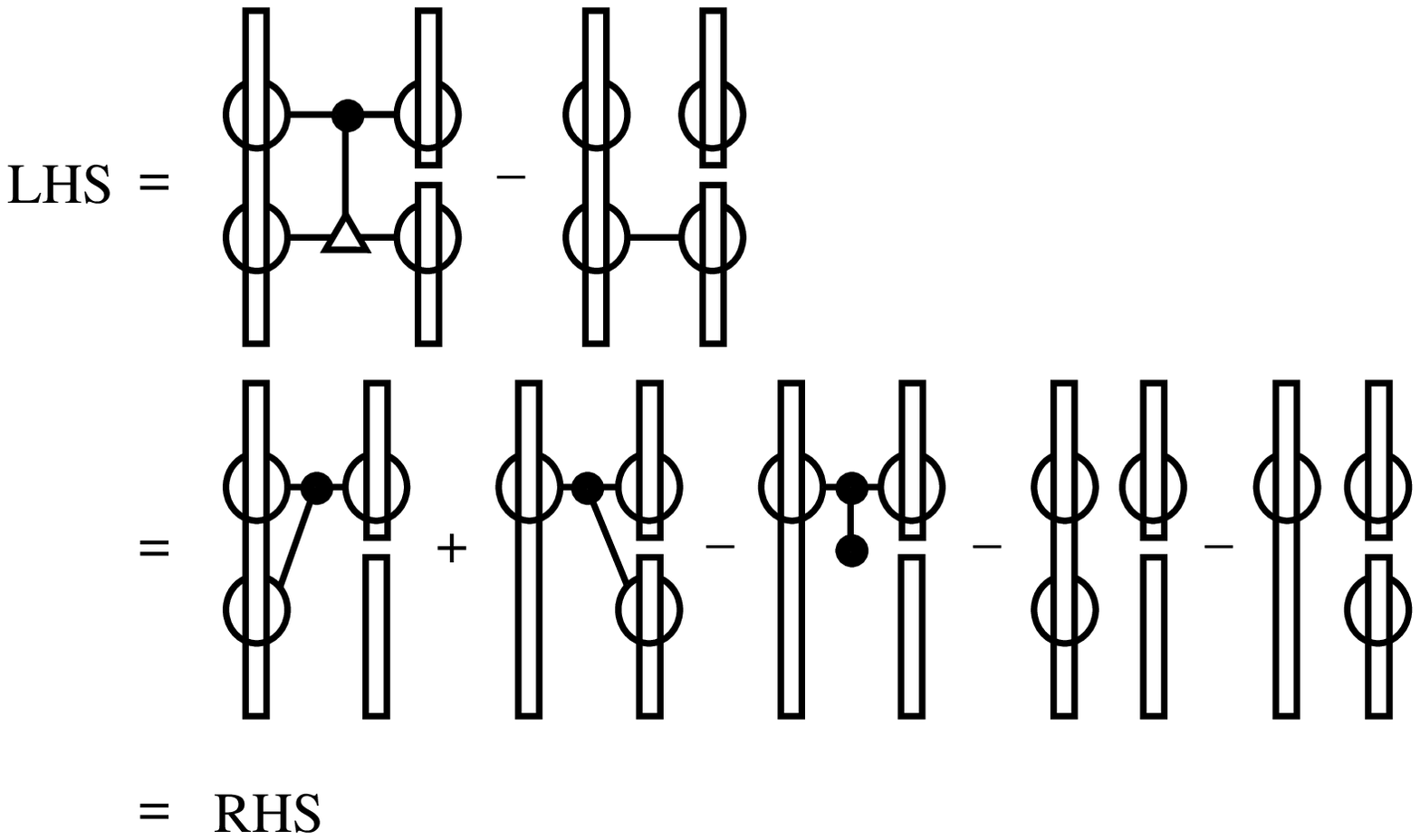}{2.8in}
$$
\caption{Proof of Lemma \ref{lemma.more.prep}}\lbl{f20}
\end{figure}

\begin{lemma}
\lbl{lemma.more} 
The relation in Figure \ref{f21} holds.
\end{lemma}

\begin{figure}[htpb]
$$ \psdraw{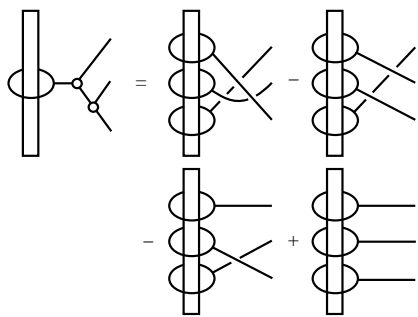}{2.2in}
$$
\caption{Breaking  two white vertices. 
An inhomogenous figure that represents
an identity in $\MM/\FF_{n+1}\MM$. The link on the left has $n$
components, and the four links on the right have $n-2$ components.}\lbl{f21}
\end{figure}

\begin{pf}
The idea of the proof is to apply Lemma \ref{lemma.preparation}
{\em twice}. Since the identity of Lemma \ref{lemma.preparation}
is inhomogenous and invloves $n$ and $n-1$ component links, whereas
the identity of Figure \ref{f21}
is inhomogenous, and invloves $n$ and $n-2$ component links,
it is not a priori clear that we can apply Lemma   \ref{lemma.preparation}
twice. Instead, we will apply the {\em proof} of Lemma  
\ref{lemma.preparation} twice. 

Now, for the proof, 
consider \lq\lq $c$'' and \lq\lq $v$'' defined as in Figure \ref{f22}.
We apply the same argument as the proof of Lemma \ref{lemma.preparation}
to the left white vertex 
in the left hand side of the required formula,
expressing the right vertex
with \lq\lq $c$'' or \lq\lq $v$'', as follows.

\begin{figure}[htpb]
$$ \psdraw{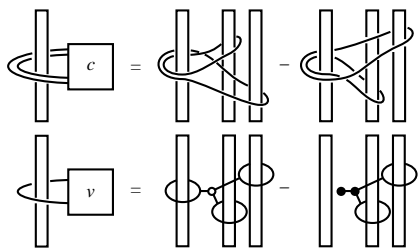}{2.2in}
$$
\caption{The definition of $c$ and $v$}\lbl{f22}
\end{figure}

Instead of $(\#1)$ in Figure \ref{f14},
we put $(\#1')$ as in Figure \ref{a1}.
We also put $(\#2')$ modifying $(\#2)$ similarly. Both $(\#1')$ and
$(\#2')$ are $n$-component links. 
In the same way as 
the former part of the proof of Lemma \ref{lemma.preparation},
we obtain $(\#1')-(\#2')$
as in Figure \ref{a2},
noting that,
in the former part,
we used, not Lemma \ref{lemma.more.prep},
but the handle slide move and calculations of alternating sum.
Note also that
we can replace $c$ with $v$
in the last picture in Figure \ref{a2}.

\begin{figure}[htpb]
$$ \psdraw{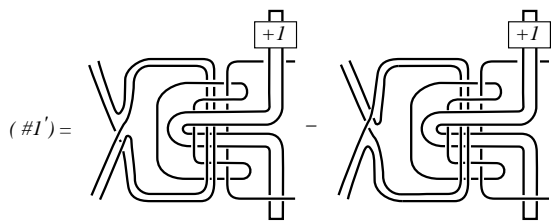}{3in}
$$
\caption{The definition of $(\#1')$}\lbl{a1}
\end{figure}

\begin{figure}[htpb]
$$ \psdraw{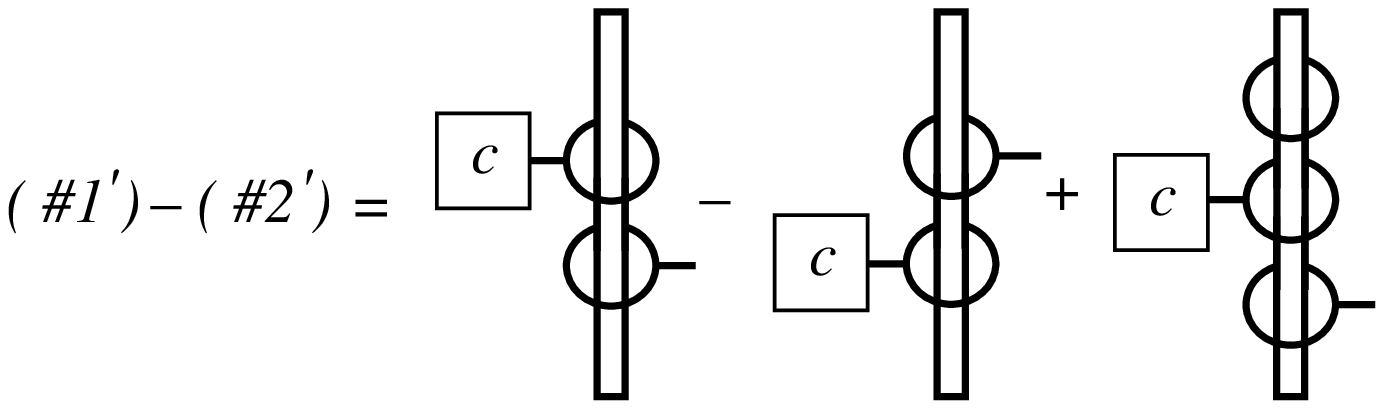}{2.5in}
$$
\caption{The calculation of $(\#1')-(\#2')$ in $\MM/\FF_{n+1}\MM$.
An inhomogenous figure where the last three links have $n-1$, $n-1$ and
$n$ components respectively.}\lbl{a2}
\end{figure}

On the other hand
we can deform $(\#1')$ as in Figure \ref{a3}
by Lemma \ref{lemma.preparation}.
Further we obtain the formula in Figure \ref{a4}
in the same way as the latter part 
of the proof of Lemma \ref{lemma.preparation}.
The same argument is valid for $(\#2')$.
Repeating the above argument once again we obtain the required formula
in the same way as the proof of Lemma \ref{lemma.preparation}.
\end{pf}

\begin{figure}[htpb]
$$ \psdraw{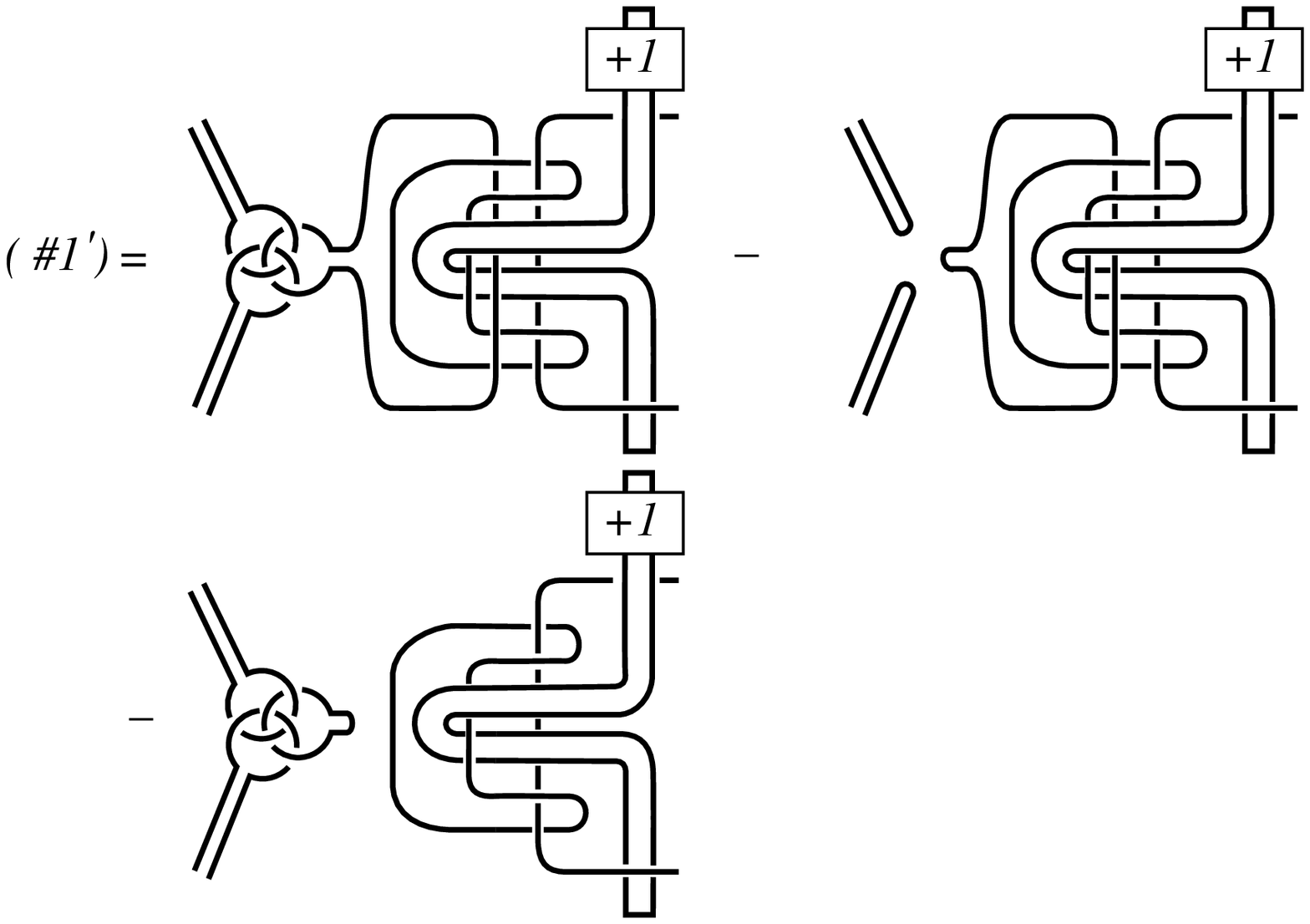}{3.2in}
$$
\caption{Another form of $(\#1')$.  A homogenous figure.}\lbl{a3}
\end{figure}

\begin{figure}[htpb]
$$ \psdraw{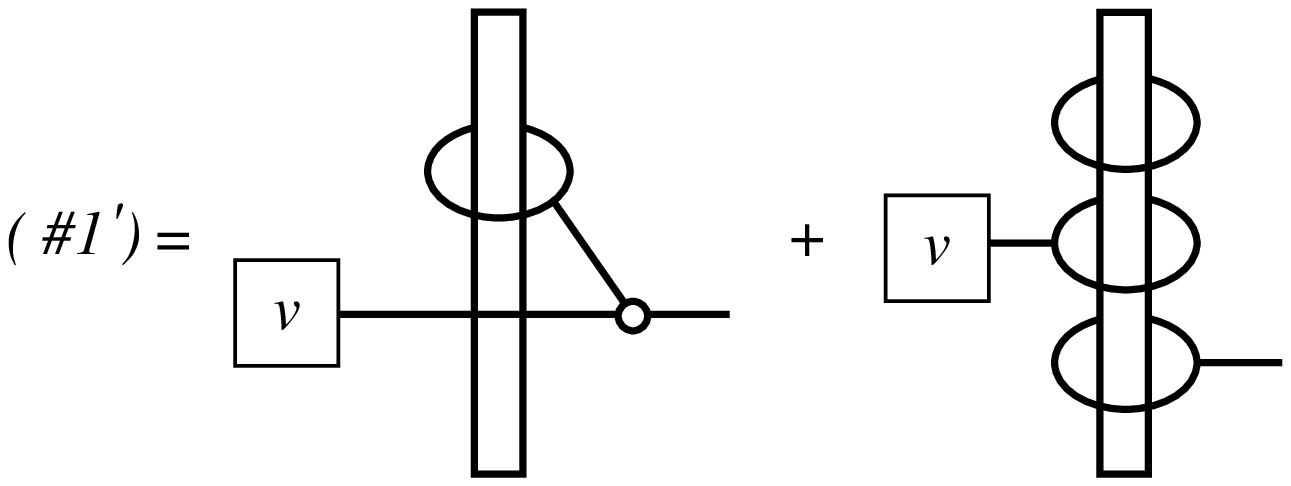}{2in}
$$
\caption{The calculation of $(\#1')$}\lbl{a4}
\end{figure}

\begin{proposition} 
The relation in Figure \ref{figure.IHX} holds,
which is called the {\em $IHX$ relation}.
\end{proposition}

\begin{figure}[htpb]
$$ \psdraw{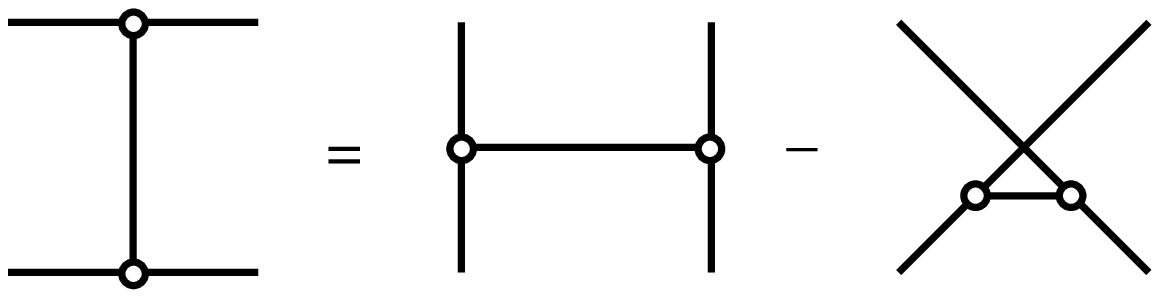}{1.5in}
$$

\caption{The $IHX$ relation}\lbl{figure.IHX}
\end{figure}

\begin{pf}
Sum up three formulas
obtained by taking cyclic permutation of both sides of the formula in 
Figure \ref{f21}.
Then we find that the right hand side vanish.
Hence we have the formula in Figure \ref{f24}, which becomes the $IHX$
 relation
using the $AS$ relation.
\end{pf}

\begin{figure}[htpb]
$$ \psdraw{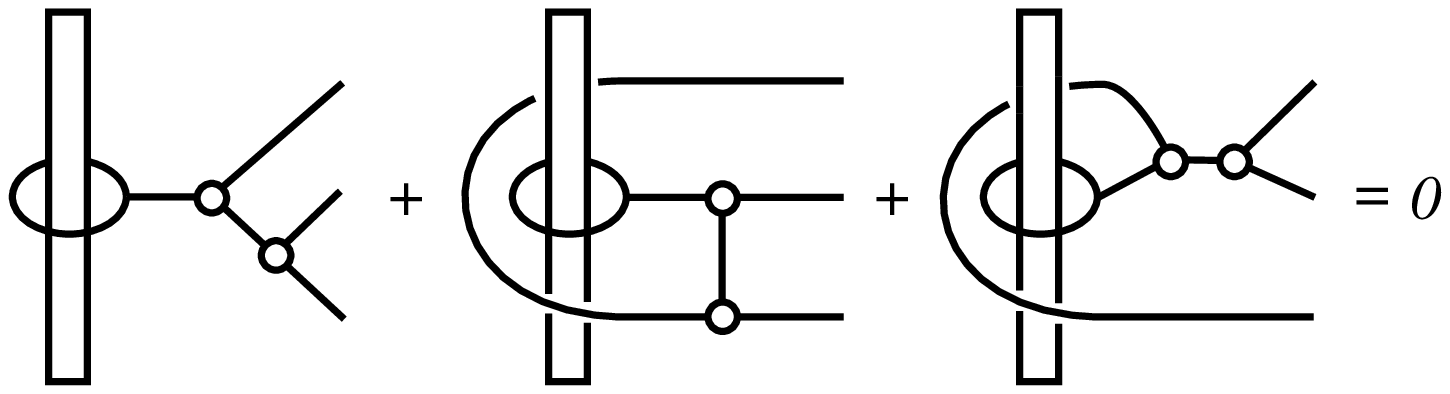}{2.2in}
$$
\caption{Proof of the $IHX$ relation}\lbl{f24}
\end{figure}

\begin{pf}[end of the proof of Theorem \ref{thm2}]
The first part of Theorem \ref{thm2} was shown above. 

For the second part, we have the fact that
$\tilde O^{\star}_m$ is onto by Theorem 1 of \OhII,
and that the deframing map $F$ is an isomorphism by its definition.
Hence it is sufficient to show that
for any ECMC $\Gamma$ with univalent and white trivalent vertices
there exists a CMC which has the same image in $\GG_m\MM$.
As in the proof of Lemma \ref{lemma.finite} in Section \ref{section.vanish},
we can show that
every connected component of $\Gamma$ is either a $Y$ component
or a trivalent graph with white vertices.
Further we can show that the image of a $\Theta$ component\footnote{
a $\Theta$ component is a trivalent graph with 2 vertices and 3 edges,
as in the Greek letter $\Theta$}
with white vertices in $\GG_m\OO$ is equal to
two times the image of a $Y$ component
by using arguments in \OhII;
note that we have the fact that $\GG_3\MM$ is one dimensional in \OhII.
Therefore we can replace a $Y$ component with half times
a $\Theta$ component with white vertices, completing this part.

For the third part, use the fact that $O_m^\star$ is onto and the definition
of $\GG_{m} \OO$. The proof of Theorem \ref{thm2} is complete.
\end{pf}


\section{Combinatorics of manifold weight systems}
\lbl{section.combi}

In this section we concentrate the combinatorics of manifold weight systems.
Our arguments are combinatorial, with little resemblance to low dimensional
topology.

We begin with the following definition:
\begin{definition}
\begin{itemize}
\item
     Let $ \cdot : \CMM \otimes \CMM \rightarrow \CMM $  be defined 
by the disjoint union of Chinese manifold characters.
\item
     Let $ \Delta : \CMM \otimes \CMM \rightarrow \CMM $  be defined
as follows: for a Chinese manifold character $\Gamma$, let
\begin{equation}
\Delta(\Gamma)= \sum\begin{Sb} \Gamma = \Gamma_1 \cup \Gamma_2 \end{Sb}
\Gamma_1 \otimes \Gamma_2
\end{equation}
where the summation is over all ways of splitting $\Gamma$ as a 
 disjoint union $\Gamma_1 \cup \Gamma_2 $, where $\Gamma_1, \Gamma_2 $
are {\em connected} Chinese manifold characters. 
\end{itemize}
\end{definition} 

\begin{pf}[of Proposition \ref{proposition.algebra}]
We claim that the above defined multiplication and commultiplication in 
$\CMM$ descends to a well defined one in $\BMM$, and that $\BMM$
becomes a commutative, co-commutative Hopf algebra.
Indeed, it follows by definition that
\begin{eqnarray*}
\Delta(AS) & = & AS \otimes 1 + 1 \otimes AS \\
\Delta(IHX) & = & IHX \otimes 1 + 1 \otimes IHX
\end{eqnarray*}
from which follows that the multiplication and the comultiplication
descend in $\BMM$. Commutativity and co-commutativity are  obvious, and so
are verifying the rest axioms of the Hopf algebra. It is an easy exercise 
to show that the primitive elements in
$\BMM$ are the {\em connected} Chinese manifold characters.
Furthermore, it follows by definition that $O_m$ is an algebra map.
The proof of Proposition \ref{proposition.algebra} is complete.
\end{pf}

Recalling that the map $\iota : \BMM \to \BBMM $ of Definition \ref{def.basic}
is an isomorphism, it follows that $\BBMM$ is also a commutative, 
co-commutative Hopf algebra, with multiplication $\tilde{\cdot}$, and
commultiplication $\tilde{\Delta}$. We can now propose the following exercise:

\begin{exercise}
Show that:
\begin{itemize}
\item
      $ \tilde{\cdot} : \BBMM \otimes \BBMM \rightarrow \BBMM $  is given
by the disjoint union of extended Chinese manifold characters.
\item
     $ \tilde{\Delta} : \BBMM \otimes \BBMM \rightarrow \BBMM $  is given
as follows: for an Chinese manifold character $\Gamma$, let
\begin{equation}
\tilde{\Delta}(\Gamma)= \sum\begin{Sb} c \in \{l,r \}^{e(\Gamma)} \end{Sb}
\Gamma_l \otimes \Gamma_r
\end{equation}
where the summation is over the set of all colorings of the edges $e(\Gamma)$
by $l,r$, and $\tilde{\Gamma}_l$ (resp. $\tilde{\Gamma}_r$) are the graphs
obtained by choosing only the $l$-colored, (resp. $r$-colored) edges of 
$\Gamma$.The graphs $\tilde{\Gamma}_l$ and $\tilde{\Gamma}_r$ have vertices
of valency $1,2$ and $3$. After splitting every vertex of valency $2$
in two vertices of valence $1$, the resulting graphs $\Gamma_l, \Gamma_r$
are Chinese manifold characters.
\end{itemize}
\end{exercise}


\ifx\undefined\bysame
	\newcommand{\bysame}{\leavevmode\hbox to3em{\hrulefill}\,}
\fi

\appendix

\section{Addendum}

After completing the text of this paper,
T.T.Q.  Le proved the following theorem,
based on our construction of the map $O_m^\star$.

\begin{theorem}\cite{Le}
\label{thm.le}
The topological invariant $\Omega(M)$
of a 3-manifold $M$ defined in \cite{LMO}
is the universal finite type invariant for integral homology
3-spheres and induces the inverse of the map \eqref{eq.Omstar}	
given in Theorem \ref{thm2}.
\end{theorem}

\begin{corollary}
\lbl{cor.cor1}
The map $O_m^\star$ of \eqref{eq.Omstar} is an isomorphism
of finite dimensional vector spaces.
\end{corollary}

Dually we have,

\begin{corollary}
\lbl{cor.cor2}
Extending \eqref{eq.O.exact}, we obtain the following short
exact sequence:
\begin{equation}\label{eq.O.exact2}
0 \rightarrow \FF_{m-1} \OO \rightarrow \FF_m \OO \rightarrow \GG_m \WW \MM
\rightarrow 0
\end{equation}
\end{corollary}

We note that in both cases of \fti s of knots and \ihs s, the existence of 
the {\em universal} \fti\ (due to Kontsevich \cite{Ko} for knots, and 
Le \cite{Le} for \ihs s)
implies the isomorphism of corollary \ref{cor.cor1} and the short
exact sequence of corollary \ref{cor.cor2}.


\begin{thebibliography}{[EMSS]}

\bibitem[AxS1]{AxS1}
S. Axelrod, I. M. Singer,  {\em Chern-Simons perturbation theory},
  Proc. XXth DGM Conference (New York, 1991) World Scientific, (1992) 3-45. 

\bibitem[AxS2]{AxS2}
S. Axelrod, I. M. Singer,  {\em Perturbative aspects of
 Chern-Simons topological quantum field theory II}, Jour. Diff. Geom. {\bf 39}
  (1994) 173-213.

\bibitem[B-N1]{B-N} D. Bar-Natan, 
        {\em On the Vassiliev knot invariants}, Topology {\bf 34} (1995)
        423-472.

\bibitem[B-N2]{B-N2} \bysame,
        {\em Computer data files} available via anonymous file transfer from
       \verb$math.harvard.edu$, user name \verb$ftp$, subdirectory \verb$dror$.
        Read the file \verb$README$ first.
   

\bibitem[Ga]{Ga} S. Garoufalidis,
        {\em On finite type 3-manifold invariants I}, M.I.T. preprint 1995,
        to appear in J. Knot Theory and its Rami.

\bibitem[GL]{GL} S. Garoufalidis, J. Levine,
       {\em On finite type 3-manifold invariants II}, Brandeis Univ. and
       M.I.T. preprint 1995.

\bibitem[Ko]{Ko} M.~Kontsevich,
	{\em Vassiliev's knot invariants},
	Adv. in Sov. Math., {\bf 16(2)} (1993), 137--150.


\bibitem[Le]{Le} T.T.Q. Le,
{\em An invariant of integral homology 3-spheres
which is universal for all finite type invariants},
preprint, 1996.

\bibitem[LMO]{LMO}T.T.Q.Le, J. Murakami, T. Ohtsuki,
{\em On a universal quantum invariant of 3-manifolds},
preprint, 1996.
 
\bibitem[Oh1]{Oh1} T. Ohtsuki,
        {\em A polynomial invariant of integral homology spheres}, Math.
        Proc. Camb. Phil. Soc. {\bf 117} (1995) 83-112.

\bibitem[Oh2]{Oh2} \bysame,
        {\em Finite type invariants of integral homology 3-spheres}, 
        preprint 1994, to appear in J. Knot Theory and its Rami.

\bibitem[Rz1]{Rz1} L. Rozansky,
   {\em The trivial connection contribution to Witten's invariant and
finite type invariants of rational homology spheres}, preprint 
{\tt q-alg/9503011}.

\bibitem[Rz2]{Rz2} L. Rozansky,
       {\em Witten's invariants of rational homology spheres at prime values
       of $K$ and the trivial connection contribution},
       preprint {\tt q-alg/9504015}.

\bibitem[Sw]{Sw} M. E. Sweedler,
        {\em Hopf Algebras}, W.A. Benjamin Inc, New York, 1969.
 
\bibitem[Va]{Va} V. A. Vassiliev,
   {\em Complements of discriminants of smooth maps}, Trans. of Math. Mono.
   {\bf 98} Amer. Math. Society., Providence, 1992.

\bibitem[Wi]{Wi}
E. Witten, {\em Quantum field theory and the Jones polynomial},
        Commun.  Math. Phys.  {\bf 121} (1989) 360-376. 


\end{thebibliography}
\end{document}